\documentclass[12pt]{article}
\usepackage{amsmath}
\usepackage{amssymb}
\usepackage{latexsym}
\usepackage{enumerate}
\usepackage[ansinew]{inputenc}
\usepackage[all]{xy}

\newtheorem{theorem}{Theorem}[section]
\newtheorem{proposition}[theorem]{Proposition}
\newtheorem{corollary}[theorem]{Corollary}
\newtheorem{lemma}[theorem]{Lemma}

\newtheorem{remark}[theorem]{Remark}

\newtheorem{problem}[theorem]{Problem}
\newtheorem{definition}[theorem]{Definition}
\newtheorem{observation}[theorem]{Observation}

\newcommand{\Z}{{\mathbb Z}}

\newcommand{\Hom}{{\rm Hom}}
\newcommand{\Aut}{{\rm Aut}}
\newcommand{\Out}{{\rm Out}}
\newcommand{\End}{{\rm End}}

\newcommand{\Ker}{{\rm Ker}\,}

\newcommand{\qed}{\quad\lower0.05cm\hbox{$\Box$}}

\def\restr{\upharpoonright}
\newcommand{\downarrowright}[1]{\downarrow
\rlap{\raise0.1cm\hbox{$\scriptstyle{#1}$}}}

\newcommand{\downarrowleft}[1]{\rlap{\kern-0.2cm
\raise0.1cm\hbox{$\scriptstyle{#1}$}}\downarrow}

\newcommand{\uparrowright}[1]{\uparrow
\rlap{\lower0.1cm\hbox{$\scriptstyle{#1}$}}}

\newcommand{\uparrowleft}[1]{\rlap{\kern-0.2cm
\lower0.1cm\hbox{$\scriptstyle{#1}$}}\uparrow}

\newcommand{\epi}{\mbox{$\to$\hspace{-0.35cm}$\to$}}

\newcommand{\mono}{\hookrightarrow}

\newcommand{\Mult}{\rm Mult}

\def\Aut{\mathop{\rm Aut}\nolimits}

\def\Out{\mathop{\rm Out}\nolimits}
\def\Inn{\mathop{\rm Inn}\nolimits}
\def\Inn{\mathop{\rm Inn}\nolimits}

\def\End{\mathop{\rm End}\nolimits}

\def\id{\mathop{\rm id}\nolimits}

\def\Hom{\mathop{\rm Hom}\nolimits}

\def\Ker{\mathop{\rm Ker}\nolimits}

\def\Z{{\mathbb Z}}

\def\fz{{\mathfrak z}}
\def\fc{{\mathfrak c}}

\def\bar{\overline }

\def\Gbar{\bar G}

\def\Hbar{\bar H}

\def\vabar{{\bar \va}}

\def\psibar{{\bar \psi}}

\def\a{\alpha}
\def\b{\beta}
\def\g{\gamma}

\def\va{\varphi}
\def\k{\kappa}

\def\n{\nu}

\def\s{\sigma}
\def\l{\lambda}

\def\aln{{\aleph_0}}
\def\al-n{{\aleph_n}}
\def\Cont{2^{\aln}}

\def\n+1d{{}^{ n+1 \downarrow }\l}

\def\size#1{\left|\,#1\,\right|}

\def\Gtilde{\widetilde G}

\def\Htilde{\widetilde H}

\def\restr{\mathop{\upharpoonright}}
\def\to{\rightarrow}
\def\arr{\longrightarrow}

\def\Pf{\smallskip\goodbreak{\sl Proof. }}
\def\Fp{\vadjust{}\penalty200 \hfill
\lower.3333ex\hbox{\vbox{\hrule\hbox{\vrule\phantom{\vrule height
6.83333pt depth 1.94444pt width 8.77777pt}\vrule}\hrule}}
\ifmmode\let\next\relax\else\let\next\par\fi \next}

\begin{document}
\setcounter{page}{1}
\title{
On localizations of quasi-simple groups with given countable center
}
\footnotetext{


The authors were partially supported by FEDER-MEC grant MTM-2016-76453-C2-1-P. Second-named author was also supported by Grant FQM-213 of the Junta de Andaluc\'{i}a.}

\author{Ram\'on Flores and Jos\'e L. Rodr\'{i}guez}
\date{\today}

\maketitle

%%%%%%%%%%%%%%%%%%%%%%%%%%%%%%%%%%%%%%%%%%%%%%%%%%%%%%%%%%%%%%%%%
%                         ABSTRACT                              %
%%%%%%%%%%%%%%%%%%%%%%%%%%%%%%%%%%%%%%%%%%%%%%%%%%%%%%%%%%%%%%%%%
\begin{abstract}
A group homomorphism $i: H \to G$ is a {\it localization} of $H$,
if for every homomorphismn $\va: H\to G$  there exists a unique endomorphism $\psi: G\to G$
 such that $i \psi=\va$ (maps are acting on the right).
 G\"obel and Trlifaj asked in \cite[Problem 30.4(4), p. 831]{GT12} which abelian groups
are centers of localizations of simple groups.
Approaching this question we show that
every countable abelian group is indeed the center of some localization
of a quasi-simple group, i.e. a central extension of a simple group.
The proof uses Obraztsov and Ol'shanskii's construction of infinite simple groups
with a special subgroup lattice and also extensions of results on
 localizations of finite simple groups by the second author and Scherer, Th\'evenaz and Viruel.
\end{abstract}

%%%%%%%%%%%%%%%%%%%%%%%%%%%%%%%%%%%%%%%%%%%%%%%%%%%%%%%%%%%%%%%%%
%%%%%%%%%%%%%%%%%%%%%%%%%%%%%%%%%%%%%%%%%%%%%%%%%%%%%%%%%%%%%%%%%
\section{Introduction}
%%%%%%%%%%%%%%%%%%%%%%%%%%%%%%%%%%%%%%%%%%%%%%%%%%%%%%%%%%%%%%%%
%%%%%%%%%%%%%%%%%%%%%%%%%%%%%%%%%%%%%%%%%%%%%%%%%%%%%%%%%%%%%%%%
Simple groups are a good source of many kinds of realization theorems.
For example, we know from Droste--Giraudet--G\"obel \cite{DGG01} that every group can be expressed as the outer automorphism of some simple group.
Using Ol'shanskii \cite{Ol91},  Obraztsov \cite{Obra96} proved that every abelian group is the center of some infinite quasi-simple group.
In \cite{Ol91} we can find Burnside groups of large  prime exponent with many extra properties. In this paper we apply Obraztsov \cite{Obra96}  to find localizations of (quasi)-simple groups with given countable center, see Theorem~\ref{RST}.

Before going into details of our results, let us recall the importance of localization, especially in Group Theory and Topology. In general, the goal of defining a localization is simplifying an object or category, sometimes by passing from global to local, sometimes by inverting morphisms.
Since the seventies, these constructions have entered with force in the categories of groups and spaces, showing their potential. Localization of spaces was introduced by Sullivan \cite{Sul74} in his proof of the Adams conjecture, Hilton--Mislin--Roitberg \cite{HMR75} introduced the group version while generalizing abelianization to nilpotent groups, Bousfield \cite{Bou75} provided the bridge between these two types of functors, and Farjoun~\cite{Far96} gave an unified treatment with the notion of $f$-localization with respect to a morphism~$f$. Note that values of group $f$-localization functors $\eta_H: H\to LH$ are localizations in the sense defined above, and conversely, given a localization $f:H\to G$, we have that  $G=L_fH$ where $L_f$ is the localization functor with respect to $f$ (see e.g.  \cite[Lemma~2.1]{Cas00}).

We recall now several subjects where these localization tools have been of great importance: concerning (co)homology of groups, the $p$-primary analysis based on $p$-localization is crucial in the modern theory of fusion systems of Broto--Levi--Oliver \cite{BLO03}, which led in particular to the surprising fact that the $\mathbb{F}_p$-homology of a group $G$ determines the homotopy structure (mod $p$) of its classifying space. Also, the exotic non-perfect group localization of a perfect group discovered in \cite{RSV05} produces a universal acyclic space \cite{BeCa99} which in turn determines the plus-construction, the key tool used by Quillen \cite{Qui73} to define algebraic $K$-theory out of the general linear groups.  In a more geometric context, localization plays a key role in the celebrated proof that all finite loop spaces are manifolds \cite{BKNP}, and very recently \cite{AIP20} localizations have been defined which are compatible with actions of groups.

More related to our approach in this paper, (co)localizations of Burnside groups (which in turn are themselves localizations of free groups) have appeared at least in two different contexts: in relation with amenability phenomena \cite{FM19}, and from the point of view of combinatorial group theory, for example in \cite{GT12} and \cite{HPR17}.
In particular, the existence of $\Cont$ varieties of groups not closed under cellular covers (another name for colocalizations) is shown in the last references, complementing the fact that there are $\Cont$ varieties closed under taking cellular covers \cite{Go}.
%; these ideas are closely to our work.
On the other hand, cellular covers of simple groups were described in \cite{BCFS13}, see also \cite{ER18}.

Our first attempt to deal with the problem of G\"obel and Trlifaj was to use Libman's work and subsequent constructions of large arbitrary single-group localizations. Libman \cite{Lib99}  showed that the natural inclusion $A_n\mono A_{n+1}$ of  alternating groups (for $n\geq 7$) is a localization. Indeed, this motivated the study of localizations between simple groups in \cite{RST02,RSV06,PS06}. First examples of infinite localizations of $A_n$ (for $n\ge 10$) were also given in \cite{Lib99}; in this case $A_n\mono SO(n)$ is a localization. In fact, for every non-abelian simple group $S$ and  cardinal $\k$, there exists a simple localization of cardinality $\ge \k$. This was proved under the (GCH) in \cite{GRS02} and extended to any universe of set theory in \cite{GS02} jointly with Shelah. Recall that  GCH (Generalized Continuum Hypothesis) states that given any infinite cardinal $\lambda$, there is no other cardinal $\lambda'$ such that $\lambda<\lambda'<2^{\lambda}$.

In this paper we will need simple groups with large Schur multipliers and thus must discard the use of these simple groups; see Karpilovsky \cite{Ka} for the restricted size of Schur multipliers of finite simple groups. Simple groups from \cite{GRS02,GS02} are here of no use, and we take advantage  of the constructions in \cite{Ol91,Obra96}. Relying on these Burnside groups, our main result Theorem \ref{main-theorem} gives a partial positive answer to Problem 4 in \cite[Chapter 30]{GT12} for countable abelian groups.

Let $F$ be a free abelian group of countable rank and $A$ be any epimorphic image. Moreover, let $H$ be a  complete, co-hopfian simple group without involutions, and suppose that $m> 10^{75}$ or $m=\infty$, such that all proper subgroups of $H$ are cyclic of an order dividing $m$, if $m<\infty$ and of infinite order otherwise. The existence of these groups $H$ follows from \cite{Obra96,Ol91}. Then we will show the following claims:
\begin{itemize}
\item[(a)] There exists a localization $H\subset G$ with $G$ complete, simple, co-hopfian $G$ and Schur multiplier $\Mult(\emph{G})\cong \emph{F}$.
\item[(b)] This localization lifts to a localization $\Htilde \to \Gtilde$ between the corresponding universal central extensions, where $\fz(\Gtilde)\cong F$.
\item[(c)] This localization also lifts to another localization $\Hbar \subset \Gbar$, between the corresponding central extensions, where $\fz{\overline G}\cong A$.
\end{itemize}

By (b) all free, abelian groups of rank $\le \Cont$ can arise as center of localization of simple groups.
More generally, we will represent countable abelian groups as centers of localizations of
(quasi)-simple groups by (c). This also leads to the following open question.

\begin{problem}
Is there a simple, complete co-hopfian group $H$ without involutions and with trivial Schur multiplier?\end{problem}

\begin{corollary}
Let $A$ be any abelian group and let $H$ be a simple, co-hopfian, complete group without involutions, of exponent  $m> 10^{75}$ or $m=\infty$, such that all proper subgroups of $H$ are cyclic of an order dividing $m$, if $m<\infty$ and of infinite order otherwise, and with Schur multiplier $\textrm{\emph{Mult}} (H)=1$.
Then there exists a localization $H \subseteq L$  with $\fz L=A$, and $L/A$ is a complete, simple and co-hopfian group of exponent $m$.
\end{corollary}

The simple groups  in Theorem \ref{main-theorem} can be investigated more thoroughly if results on finite simple groups obtained in  \cite{RST02} and \cite{RSV06} are first extended to infinite groups. It is important to observe that localizations $H\to G$ between two simple groups often induce localization $\Aut(H) \to \Aut(G)$, as will be discussed next. We will also consider the induced homomorphism $\widetilde H \to \widetilde G$ between their corresponding universal and intermediate central extensions.

\noindent
{\bf Acknowledgments:} This paper is dedicated to the memory of professor R\"udiger G\"obel.

\setcounter{equation}{0}

%%%%%%%%%%%%%%%%%%%%%%%%%%%%%%%%%%%%%%%%%%%%%%%%%%%%%%%%%%%%%%%%%
%%%%%%%%%%%%%%%%%%%%%%%%%%%%%%%%%%%%%%%%%%%%%%%%%%%%%%%%%%%%%%%%%
\section{Localization of infinite simple groups}
%%%%%%%%%%%%%%%%%%%%%%%%%%%%%%%%%%%%%%%%%%%%%%%%%%%%%%%%%%%%%%%%
%%%%%%%%%%%%%%%%%%%%%%%%%%%%%%%%%%%%%%%%%%%%%%%%%%%%%%%%%%%%%%%%
First we extend some results from \cite{RST02} to infinite simple groups.   If $G$ is a group and $g,h\in G$, then $hg^*= h^g=g^{-1}hg$ is conjugation and $c:G\arr \Aut G \ (g\mapsto g^*)$ is the conjugation map with kernel $\fz G$, the center of $G$. As usual, let $G^*=\{g^*\mid g\in G\}=\Inn G\subseteq \Aut G$. It is immediate for $\a\in \Aut G$ that the following diagram commutes.
%\begin{figure}[placement specifier]
%\label{diagram}
$$\xymatrix{
 G \ar[d]_{\cong}^\alpha \ar[r]^c & \Aut(G) \ar[d]_\cong^{\alpha^*} \\
  G \ar[r]^c           & \Aut(G).     }
$$
%\caption{Diagram for (c) and (d);}
%\end{figure}
That is, for any $g\in G$ we have
\begin{eqnarray} \label{conjugates} g^*\a^*= (g\a)^* \text{ or equivalently } (gc)\a^*=(g\a)c, \end{eqnarray}

because $gc\a^*=g^*\a^*= \a^{-1}g^*\a$, $(g\a)c = (g\a)^*$ and for any $h\a\in G$ we have $$(h\a)\a^{-1}g^*\a= h g^*\a = (g^{-1}hg)\a= (g\a)^{-1}(h\a)(g\a)= (h\a)^{g\a}=(h\a)(g\a)^*.$$

Now suppose that $G,H$ are groups with trivial center.

 Let $H\subseteq G$, and $h\in H$. In order to emphasize \emph{conjugation by $h$ extended to $G$}, we also write $h^{*G}\in G^*$, thus $h^*= h^{*H}\subseteq h^{*G}$.
\begin{definition}\label{ext} If $j: \Aut H\arr \Aut G$, we say that \emph{$j$ extends  $H^*\subseteq G^*$}, if $h^*j=h^{*G}$ for all $h\in H$.\end{definition}

An easy argument from \cite[Lemma 1.2]{RST02} can be used for a more general
\begin{lemma} \label{liftaut} Suppose that $G,H$ are groups with trivial center. If $j: \Aut H\arr \Aut G$ extends $H^*\subseteq G^*$, then also any $\a\in\Aut H$ extends to $\a j\in \Aut G$, i.e. $(\a j)\restr H=\a$.\end{lemma}

\Pf If $h\in H$, then we must show that $h \a =h(\a j)$, that is, the left square in the following cube commutes.
$$
\xymatrix{
& H \ar[dl]_\alpha \ar[rr]^c \ar@{->}'[d][dd]_i & & \Aut(H) \ar@{->}[dd]^j \ar@{->}[dl]_{\alpha^*}
\\
H \ar[rr]^c \ar[dd]_i & & \Aut(H) \ar[dd]_(.6)j &
\\
 & G \ar[dl]_{\alpha j}  \ar@{->}'[r][rr]^c&  & \Aut(G) \ar[dl]_{(\alpha j)^*}
 \\
G \ar[rr]^c & & \Aut(G) &
}
$$
First we apply (\ref{conjugates}) twice to $H$ and $G$, respectively, and get $(h\a)^* = h^*\a^*$ as well as $(h(\a j))^{*G} = h^{*G}(\a j)^*$ (i.e. top and bottom squares commute). Moreover, since $j$ is a homomorphism, it follows $(h^*\a^*)j = h^{*G}(\a j)^*$ (i.e. right square commutes). Recall now that $\a\in \Aut H$. From these three equalities and the assumption $(h\a)^*j = (h\a)^{*G}$ on $j$ we get
$(h(\a j))^{*G} = h^{*G}(\a j)^*=(h^*\a^*)j = (h\a)^*j = (h\a)^{*G},$
hence $(h\a)^{*G}= (h(\a j))^{*G}$. Since $G$ has trivial center, it follows $h\a= h(\a j)$ and $(\a j)\restr H = \a$.
%$$
%\xymatrixcolsep{5pc}
%\xymatrix{
%H \ar@{_{(}->}[dl]_\alpha \ar[rr]\ar@{->}'[d][dd]_i & & Aut(H) \ar@{->}[dd]^j \ar@{->}[dl]_{\alpha^*}
%\\
%H  \ar@{->}[rr] \ar@{->}[dd]_i &  & Aut(H)  \ar@{->}[dd]^j
%\\
%& G  \ar@{->}[dl]_{\alpha j} & & \Aut G   \ar@{.>}[dl]_(.6)^{(\alpha j)^*
%\\
%G  \ar@{->}[rr]^(.6)^c &  &  \Aut G
%}
%$$
\Fp

\bigskip

The lemma will be used to characterize certain localizations $H\subseteq G$. When $G,H$ are finite simple groups, then the next Theorem \ref{RST} is due to \cite[Theorem 1.4]{RST02} (and a variation of \cite[Corollary 4]{GRS02}), because finite simple groups are obviously co-hopfian. Recall that a group $G$ is \emph{co-hopfian}, if every monomorphism $G\hookrightarrow G$  is an automorphism of $G$. In one direction of our extension we will need that $G$ is co-hopfian.

\begin{theorem} \label{RST} Let  $H\subseteq G$ be an extension of two simple groups, and assume that $G$ is co-hopfian. Then, $H\subseteq G$ is a localization, if and only  if the following conditions hold.
\begin{enumerate}
\item[(1)]  There is a monomorphism $j:\Aut H \mono \Aut G$ which extends $H^*\subseteq G^*$.
\item[(2)] Any subgroup of $G$ isomorphic to $H$ is conjugate to $H$ in $\Aut G$.
\item[(3)] The centralizer $\fc_{\Aut G} H=1$ is trivial.
\end{enumerate}
%Conversel, if the extension $H\subseteq G$ of simple groups satisfies $(1)-(3)$ and $G$ is co-hopfian, then %$H\subseteq G$ is a localization.
\end{theorem}
\Pf  Assume that $H\subseteq G$ is a localization. Denote the inclusion by $i: H \mono G$.

(1) If $\a\in \Aut H$, then for $\a i: H\arr G$ there is a unique $\b\in \End G$ such that $\a = \b\restr H$. Since $G$ is simple, the map $\b$ is injective. The map $\a^{-1}i: H\arr G$ also has a unique extension $\g\in \End G$ such that $\a^{-1}= \g\restr H$, and from $\a^{-1}\a=\id_H$ and uniqueness of localizations it follows $\g\b=\id_G$. Hence $\b$ is also surjective and $\b\in \Aut G$. It is clear by the uniqueness of localizations that the map
$$j:\Aut H\arr \Aut G \ (\a\mapsto \a j =\b)$$
is a monomorphism.  Note that $h\in H$ by definition of $j$ is mapped to $h^*j= h^{*G}$, hence $j$ extends $H^*\subseteq G^*$ and $(1)$ holds.

(2) If $K$ is a subgroup of $G$ and isomorphic to $H$, then there is some monomorphism $\a:H\arr K\subseteq G$. By the localization there is $\b\in \End G$ such that $\b\restr H=\a$, and since $G$ is co-hopfian, it follows that $\b\in \Aut G$. Hence $\b^{-1}K\b= K^{\b^*}=H$, as required.

(3) If $\b\in \Aut G$, and $\a=\b\restr H=\id_H$, then by the uniqueness of extensions it follows $\b=\id_G$, hence (3) holds.

Conversely, we assume that $H,G$ are simple, $G$ is also co-hopfian and we want to show that $H\subseteq G$ is a localization, if $(1),(2),(3)$ hold.  If $\va\in \Hom(H,G)$, then we must show that there is a unique $\b\in \End G$ such that $\va= \b\restr H$. If $\va$ is the trivial homomorphism, then we choose $\b$ also trivial and note that $\b$ is unique, because $G$ is simple and $0\ne H\subseteq \Ker \b$ forces $G\b$ to be trivial too. We may assume that $\va\ne 0$, and $\va$ is injective, because also $H$ is simple. If $K=H\va$, then there is  an isomorphism $\a:H\arr K$, and by $(1)$ and $(2)$ we find $\b\in \Aut G$ such that $K^{\b^*}=H$. The automorphism $\a\b\in \Aut H$ extends to an automorphism $\psi\in \Aut G$ by (1) (using Lemma \ref{liftaut}). Then $\psi\b^{-1}$ is the desired extension of $\a $.
Its uniqueness follows from $(3)$ and the assumption that $G$ is co-hopfian. \Fp \bigskip

We will apply the following immediate consequence of Theorem \ref{RST}.

\begin{corollary}
\label{Criterion-complete}
Let $H\subseteq G$ be an extension of simple, complete groups.  If $H\subseteq G$ is a localization,
then the following two conditions hold.
\begin{enumerate}
\item[(1)] Any subgroup of $G$ which is isomorphic to $H$ is conjugate to $H$.
\item[(2)] $\fc_G(H)=1$.
\end{enumerate}
Conversely, if $G$ is co-hopfian and $(1)$ and $(2)$ hold, then $H\subseteq G$ is a localization.
\end{corollary}

%%%%%%%%%%%%%%%%%%%%%%%%%%%%%%%%%%%%%%%%%%%%%%%%%%%%%%%%%%%%%%%%%
%%%%%%%%%%%%%%%%%%%%%%%%%%%%%%%%%%%%%%%%%%%%%%%%%%%%%%%%%%%%%%%%%
\section{Automorphism groups and localizations}
%%%%%%%%%%%%%%%%%%%%%%%%%%%%%%%%%%%%%%%%%%%%%%%%%%%%%%%%%%%%%%%%%
%%%%%%%%%%%%%%%%%%%%%%%%%%%%%%%%%%%%%%%%%%%%%%%%%%%%%%%%%%%%%%%%%
We next extend some results from Section 2 of \cite{RSV06}
to simple groups of arbitrary size.
Recall that $\Out G$ denotes the outer automorphism group $\Aut G/G$ of a simple group $G$ (where we identify $G$ and $G^*$).

Lemma 2.1 of \cite{RSV06} holds for infinite groups, too.
\begin{lemma}
\label{subgroups-aut}
Let $G$ be a non-abelian simple group. Then any proper normal subgroup of $\Aut G  $ contains  $G$.
In particular, any endomorphism of $\Aut G$ is either a monomorphism or contains $G$ in its kernel. \Fp
\end{lemma}

Lemma 2.2 of \cite{RSV06} states that any non-abelian finite simple subgroup of $\Aut G$ is contained in $G$. The easy argument is based on the solution of the Schreier conjecture, which ensures that $\Out G$ is solvable for every finite non-abelian simple group $G$.
But for infinite groups this is not longer true, because from
\cite{DGG01} we know that all groups are outer automorphism groups of simple groups.
In order to proceed to infinite groups we will assume that the outer automorphism group is hyperabelian, which extends solvability. Recall that a group $A$ is called {\it hyperabelian} if every non-trivial epimorphic image of $A$ has a non-trivial abelian normal subgroup. In particular, $A$ can be filtered as an ascending (possibly transfinite) union $\bigcup K_i$ of normal subgroups with abelian factors $K_{i+1}/K_i$.
Hence, $\Hom(H,A)=0$ for all non-abelian simple groups $H$.

\begin{lemma}
\label{simple-subgroup-Aut}
If $G$ is a simple group and $\Out G$ is hyperabelian, then any non-abelian simple subgroup of $\Aut G$ is contained in $G$.
\end{lemma}
\Pf
We identify $G$ and $G^*$, which is normal in $\Aut G$. If  $H\neq 1$ is a simple non-abelian subgroup of $\Aut G$,
then $H\cap G$ is normal in $H$. Hence it is either $1$ or $H$ by simplicity of $H$. If $H\cap G=1$ then $H\cong GH / G$ is a subgroup of $\Out G$ which, by assumption, is hyperabelian, hence $H$ is cyclic, a contradiction. Otherwise,  $H\cap G=H$, hence $H\subseteq G$ as required. \Fp \bigskip

Next we want to extend \cite[Theorem 2.3]{RSV06}. We will follow the arguments in \cite[p. 770]{RSV06}, but must adjust two of the conditions. Since $\Aut G$ need not be complete, we must replace $\Aut G$ by  $\Aut(\Aut G)$. Moreover, (similar to the above,) we assume that $\Aut G$ is co-hopfian and will apply Lemma \ref{simple-subgroup-Aut}.

\begin{theorem}
\label{criterium-Aut-loc}
Let $j:\Aut H \mono \Aut G$ be an inclusion of the automorphism groups of two non-abelian simple groups $H$ and $G$.
Assume that $\Aut G$ is co-hopfian and $\Out G$ is hyperabelian. Then $j:\Aut H \mono \Aut G$ is a localization, if and only if the following four conditions are satisfied.
\begin{itemize}
\item[(a)] If $\Omega =\{k: \Aut H \mono  \Aut G\}$ denotes all monomorphisms of $\Aut H$ in $\Aut G$, then $\Aut(\Aut G)$ acts transitively on $\Omega$.
    \item[(b)] $\fc_{\Aut(\Aut G)}\Aut H   =1$.
\item[(c)] Any homomorphism $\va': \Aut H / (\Aut H\cap G) \to \Aut G$ extends uniquely to a $\psi': \Out G \to \Aut G$ such that $ k \psi' = \pi \varphi'$, where $\pi: \Out H \epi \Aut G  / \Aut H \cap G$ is  the canonical projection (see Figure 1).
\item[(d)] If a homomorphism $\va: \Aut H\to \Aut G$ contains $H$ in its kernel, then also $\Aut H\cap G \subseteq \Ker \va $.
\end{itemize}
\end{theorem}
\Pf Suppose that $j: \Aut H\mono \Aut G $ is a localization and let $\va: \Aut H  \mono \Aut G  $ be any monomorphism. Then $\va$ extends to a homomorphism $\psi: \Aut G   \to \Aut G$ such that $j \psi  = \va$. We must show that $\psi \in \Aut G  $.

If $\Ker \psi \neq 0$, then $G\subseteq \Ker\psi$ by Lemma \ref{subgroups-aut}, and thus $H\subseteq \Ker\psi$, which is impossible since $\va$ is monomorphism. Thus  $\psi$ is injective and also surjective, because $\Aut G  $ is co-hopfian. It follows that $\psi$ is an automorphism.
Hence (a) and (b) follow clearly.

To show (c), let $\va': \Aut H / \Aut H \cap G \to \Aut G  $ be any homomorphism (see Figure~1). Let $\va=p \va' : \Aut H  \epi \Aut H / \Aut H \cap G \to \Aut G  $.
Since $j$ is a localization, there exists a unique $\psi: \Aut G   \to \Aut G  $ such that
$j \psi = \va$. And since $\va$ is not injective, also $\psi$ is not injective, and $G\subseteq \Ker\psi$  by Lemma \ref{subgroups-aut}.
Hence $\psi$ factors through $\Out G $ and there is a unique homomorphism
$\psi': \Out G  \to \Aut G  $ such that $p \psi' = \psi$.

We claim $\pi\va'=k\psi'$.

Composing the left-hand side with $p$ we have $p\pi\va' =p \va'= \va$ and on the right-hand side by the above $p k\psi'=j p\psi'=j\psi=\va$, hence $p\pi\va' =p k\psi'$. Since $p$ is surjective, it follows  $\pi \va' =k\psi' $.

Suppose that there exists another $\psi'':\Out G \to \Aut G$ such that $k\psi'' = \pi \va'$.
Then $p k \psi'= p k \psi''$, and therefore $j p\psi' = j p\psi''$. Since $j$ is a localization, it follows $p\psi'=p \psi''$, and therefore $\psi'=\psi''$, because $p$ is an epimorphism.
Hence (c) holds.

\begin{figure}
\label{diagram}
$$
\xymatrixcolsep{5pc}\xymatrix@!0{
& H \ar@{_{(}->}[dl] \ar@{^{(}->}[rr]\ar@{->}'[d][dd] & & G \ar@{->}[dd] \ar@{->}[dl]
\\
\Aut H \cap G \ar@{->}[rr] \ar@{->}[dd] &  & G  \ar@{->}[dd]
\\
& \Aut H  \ar@{=}[dl] \ar@{->}[dr]^\va \ar@{^{(}->}'[r][rr]^j \ar@{->>}'[d][dd]^p & & \Aut G   \ar@{->>}[dd]^p \ar@{.>}[dl]_(.6)\psi
\\
\Aut H  \ar@{->}[rr]^(.6)\va \ar@{->>}[dd]_p&  &  \Aut G   \ar@{->>}[dd]^>p
\\
&\Out H    \ar@{->>}[dl]_(.5)\pi \ar@{->}'[r][rr]^(.3)k && \Out G  \ar@{->}[dl] \ar@{.>}[ul]_{\psi'}&
\\
\Aut H / \Aut H \cap G \ar@{->}[rr] \ar@/_/[uurr]_{\va'} && \Out G
}
$$
\caption{Diagram for (c) and (d).}
\end{figure}

For (d) suppose that $\va: \Aut H  \to \Aut G  $ is a homomorphism with $H\subseteq \Ker \va$ and (by the localization $j$) choose $\psi: \Aut G   \to \Aut G$ such that  $j \psi = \va $. Then, $\Ker \psi \supseteq (\Ker \va)j \neq 0$ and, by Lemma \ref{subgroups-aut}, $G\subseteq \Ker \psi$, and therefore, $\Aut (H)\cap G\subseteq \Ker \va$.

Suppose conversely that $(a)-(d)$ hold, and let $\va: \Aut H  \to \Aut G  $ be any group homomorphism.
If $\Ker(\va)=0$, then $\va \in \Omega$ and by (a) there exists an automorphism
$\psi: \Aut G  \to \Aut G  $ such that  $j \psi = \va$. Uniqueness is guaranteed by (b).

If $\Ker(\va)\neq 0$, then $H\subseteq \Ker(\va)$  by Lemma \ref{subgroups-aut}.
And by (d) also $\Aut H  \cap G\subseteq \Ker(\va)$, and therefore, $\va$ factors
through some homomorphism $\va': \Aut H / \Aut H \cap G \to \Aut G  $ such that $p\va'  =p\pi\va' = \va$. By (c) this homomorphism
extends uniquely to some $\psi': \Out G  \to \Aut G  $ such that $k\psi'  = \pi\va'$.

 We next show that the composition $\psi:=p\psi' $ satisfies $j\psi = \va $. Indeed,
$$j p \psi'=p k \psi' = p \pi \va' =\va.$$

For the uniqueness suppose that another homomorphism $\widehat\psi: \Aut G   \to \Aut G  $
satisfies $ j\widehat\psi = \va$. Then  $\Ker\widehat\psi\neq 0$ and $G\subseteq \Ker(\widehat\psi)$
by Lemma \ref{subgroups-aut}, hence $\widehat\psi$ factors through some $\widehat\psi': \Out G  \to \Aut G  $
 such that $p \widehat\psi' =\widehat\psi$.
As $p: \Aut H \epi \Out H  $ is an epimorphism, and
$$p k \widehat\psi'= j p\widehat\psi'  = j \widehat\psi  = \va=  p \pi \va'$$
we have $k \widehat\psi' = \pi\va'$. By uniqueness of (c) we get $\widehat\psi' = \psi'$,
and therefore $\widehat\psi=\psi$, as desired. \Fp \bigskip

It is clear that the following holds.

\begin{observation}
\label{localization-into-cohopfian}
If $H\subseteq G$ is a localization of non-abelian simple groups, and $G$ is co-hopfian, then every monomorphism
$\va: H \mono G$ extends to a unique $\psi\in \Aut G$ such that $\psi\restr H = \va$.
\Fp
\end{observation}

If the assumption  of co-hopfian is removed, localizations $\va: H \mono G$  for which the previous observation does not hold can be found, as it is shown in the following (abelian) examples: the localization $\Z\mono \Z[1/2]$, where $\Z\stackrel{3}{\mono} \Z[1/2]$, does not extend to any automorphism,
and  neither does the isomorphism
$\Z[1/2] \stackrel{2}{\to}\Z[1/2]$ restrict to an isomorphism on $\Z$.

A version of  \cite[Theorem 2.4]{RSV06} also extends to infinite groups (by adding co-hopfian).

\begin{theorem}
Let $H\subseteq G$ be a localization of non-abelian simple groups, suppose that $G$ and $\Aut G  $ are co-hopfian, $\Aut H $ is complete and $\Out G $ is hyperabelian.
Then any extension $j:\Aut H \mono \Aut G  $ is a localization,
if and only if conditions (b), (c) and (d) as in Theorem \ref{criterium-Aut-loc} are satisfied.
If also $\Aut G  $ is complete, then $j$ is a localization, if and only if (c) and (d) hold.
\end{theorem}
\Pf
We apply the previous Theorem \ref{criterium-Aut-loc}
for the converse implication.
Let $\va: \Aut H  \mono \Aut G  $ be any monomorphism decomposed as
$$\Aut H  \stackrel{\va}{\to} \va(\Aut H ) \stackrel{k}{\mono}\Aut G  ,$$
and let $\va': H \stackrel{\cong }{\to} H' \stackrel{k}{\mono} G$ be its restriction.
If the localization $H\subseteq G$ is denoted by $i$, then by Observation \ref{localization-into-cohopfian}, there exists
a unique  $\alpha\in \Aut G$ such that $\va' k\alpha  = i$. And clearly, conjugation $\alpha^*: \Aut G   \to \Aut G  $
restricts to an automorphism $\alpha^*: \Aut(H')\to \Aut H $ such that $k \alpha^*  =  \alpha^* j$.
By completeness the composite $\va \alpha^*: \Aut H  \to \Aut H $ is conjugation $\beta^*$
by some automorphism $\beta\in \Aut H $.
Now $ \beta^* j= ((\beta)j)^* j$, by uniqueness of localization, and then
$$
j ((\beta)j)^* (\alpha^*)^{-1} =  \beta^* j (\alpha^*)^{-1}   = \beta^* (\alpha^*)^{-1} k =  \va,
$$
therefore $\psi=((\beta)j)^*(\alpha^*)^{-1}$ extends $\va$, i.e. $j \psi = \va$, as desired.

Now suppose that $\Aut G  $ is complete. If $\beta^*: \Aut G  \cong \Aut G  $ is an inner automorphism
that centralizes $\Aut H $, then it centralizes $H$, and therefore $\beta: G\to G$ is
the identity homomorphism, by uniqueness of the localization $H\subseteq G$. \Fp \bigskip

%%%%%%%%%%%%%%%%%%%%%%%%%%%%%%%%%%%%%%%%%%%%%%%%%%%%%%%%%%%%%%%%%
%%%%%%%%%%%%%%%%%%%%%%%%%%%%%%%%%%%%%%%%%%%%%%%%%%%%%%%%%%%%%%%%%
\section{Universal central extensions and localizations}
%%%%%%%%%%%%%%%%%%%%%%%%%%%%%%%%%%%%%%%%%%%%%%%%%%%%%%%%%%%%%%%%
%%%%%%%%%%%%%%%%%%%%%%%%%%%%%%%%%%%%%%%%%%%%%%%%%%%%%%%%%%%%%%%%
We next extend results from Section 1 of \cite{RSV06}
to non-abelian simple groups of arbitrary size.

Corollary \ref{localization-superperfect} below was stated in \cite{RSV06}
only for finite groups and was used to find many examples of non-simple localizations of finite non-abelian simple groups, thus showing that simplicity
is not preserved under localization functors in general.
For infinite groups this corollary could be useful to realize
any abelian group as the center of some simple localization (see
Corollary \ref{localization-superperfect} and Question  \ref{Ol}).

Let $G$ be a non-abelian simple group and denote by $\textrm{Mult}({G})=H_2(G;\Z)$
its Schur multiplier. Consider
$$0\to \textrm{Mult}(G) \hookrightarrow \Gtilde \epi G \to 1$$
the universal central extension of $G$.
Recall that the group $\Mult(\emph{G})$ coincides with the center $\fz \Gtilde$ of $\Gtilde$.

If $G$ is a non-abelian simple group, then any proper normal subgroup of $\Gtilde$
is contained in $\Mult(\emph{G})$. In particular, any non-trivial endomorphism of $\Gtilde$ is either a monomorphism,
or its kernel is contained in $\Mult(\emph{G})$.

\begin{proposition}
Let $G$ be any group and let $H \subset G$ be a maximal subgroup.
Suppose that $H$ is non-abelian simple and $G$ has some element of order $p$, where $p$ is the index of $H$ in $G$.
If the inclusion $H\hookrightarrow G$ is a localization, then $G$ is simple.
\end{proposition}
\Pf For the convenience of the reader we repeat the proof of \cite[Proposition 1.1]{RSV06}.
Suppose that $H$ is normal in $G$. Since $H$ is maximal, $G/H\cong C_p$. By
hypothesis, there exists an element $x\in G$ of order $p$, yielding a non-trivial endomorphism
$\psi: G\to G/H \cong C_p \mono G$. But $H\psi=0$, hence $\psi=0$, because $H\mono G$
is a localization, which is a contradiction. So $H$ is not normal in $G$ and
suppose, for contradiction, that $G$ has a proper normal subgroup $N$.
Then $N\cap H=1$, otherwise $N\cap H = H$ would imply $N=G$ by the maximality of $H$
and the fact that $H$ is not normal. Therefore, $G/N\cong H$.
But then the identity homomorphism $\id: G\to G$ and the projection $G\epi G/N$ extend
an inclusion $H\mono G$, which contradicts the uniqueness of localizations. \Fp \bigskip

The following will be applied in our main Theorem \ref{main-theorem}.
\begin{proposition}
\label{mult-mult}
Let $i:H\mono G$ be an inclusion of two non-abelian simple groups, and
$j:\Htilde  \to \Gtilde$ be the induced homomorphism
between the universal central extensions.
%Moreover, let $\overline G \epi G$ be any intermediate central extension,
%$\overline H\epi H$ be the induced cover along $i:H\mono G$ (taking pull-back), and
%$k:\overline H \mono \overline G$ be the induced homomorphism.
If every homomorphism  $\Htilde  \to \Gtilde$ sends
$\textrm{\emph{Mult}}(H)$ into $\textrm{\emph{Mult}}(G)$, then $i$ is a localization, if and only if $j$ is a localization.
\end{proposition}
\Pf  As in \cite[Proposition 1.2]{RSV06},
we show first that $p:\Gtilde \epi G$ and $q:\Htilde \epi  H$ induce
a bijection $F:\Hom(\Htilde ,\Gtilde )\longrightarrow \Hom(H,G)$ characterized
by the property that $F(\tilde \va)$ (with $\tilde \va \in \Hom(\Htilde ,\Gtilde ))$ is the unique morphism $\va:H\to G$
such that $ \tilde\va p=q \va$.
If $\tilde\va : \Gtilde  \to \Htilde $ is any homomorphism, then by hypothesis $(\Mult(\emph{H}))\tilde\va\subseteq \Mult(\emph{G})$,
hence it induces a unique homomorphism $\va: H\to G$ such that $\tilde\va p= q \va$.
To show that $F$ is surjective, we use that the universal central extension
is a cellular cover (i.e. a co-localization functor), hence any homomorphism
$\va: H\to G$ induces a unique homomorphism $\tilde \va: \Htilde  \to \Gtilde$
such that $ \tilde\va p=q \va $. \Fp \bigskip

The next results are parallel to \cite[Corollary 1.7, Proposition 1.4]{RSV06}.
\begin{corollary}
\label{trivial-schur}
Let $i :H \mono G$ be an inclusion of two non-abelian simple groups. Moreover, let
H be superperfect, $j :H =\Htilde \mono G$ be the induced homomorphism on the
universal central extensions. Then $i :H \mono G$ is a localization if and only if $j :H \mono \Gtilde$ is a localization.
\Fp
\end{corollary}

The universal property of the universal central extension yields the following well-known fact:
\begin{corollary}
If $G$ is a non-abelian simple group and $p:\tilde G\epi G$ is its universal central extension, then $ \Aut(\tilde G) \cong \Aut G  $.
\end{corollary}

\begin{proposition}
\label{universalcover-localization}
Let $G$ be a non-abelian simple group with non-trivial Schur multiplier. Then, $p: \Gtilde \epi G$ is a localization,
if and only if there is no monomorphisms $\Gtilde/N \mono G$ for any normal subgroup $N\varsubsetneq \textrm{\emph{Mult}}(G)$ of $\Gtilde$.
\end{proposition}
\Pf
Suppose that $\Gtilde \epi G$ is a localization and
suppose that there exists a monomorphism $\va: \Gtilde/N \mono G$ for some $N\neq \Mult(\emph{G})$.
By the localization property there is some endomorphism
$\va':G\to G$ such that $p \va' =q \va$, where $q: \Gtilde \epi \Gtilde/N$ is the canonical projection.
But then, for $x \in \Mult(\emph{G})\setminus \emph{N}$, we get $1 \neq x q \varphi = x p \varphi' =1$, which is a contradiction.

%$(\Mult(\emph{G}))\va=(\Mult(\emph{G}))p\va'=1$, contradicting that $\va$ is a monomorphism.

Conversely, let $\va:\Gtilde \to G$ be any homomorphism.
If $\va=0$, then obviously $\va$ extends to $\va'=0$ on $G$.
If $\va\neq 0$, then $\Ker \va \subseteq \Mult(\emph{G})$, and by hypothesis $\Ker\va= \Mult(\emph{G})$
and hence $\va$ factors through a homomorphism $\va': G\to G$.
Uniqueness of $\va'$ is guaranteed in any case because $p: \Gtilde \epi G$ is an epimorphism.
\Fp

\bigskip

Finally \cite[Theorem 1.5]{RSV06} extends to infinite simple groups.

\begin{theorem}
\label{non-split}
Let $i:H\mono G$ be an inclusion of two non-abelian simple groups and
$j:\Htilde \to \Gtilde$ be the induced
homomorphism on the universal central extensions. Assume that  $G$ does not contain any non-split
central extension of $H$ as a subgroup. Then $i:H\mono G$ is a localization,
if and only if $j:\tilde H \to \tilde G$ is a localization.
\end{theorem}
\Pf
Let $\varphi: \Htilde \to \Gtilde$, then
$\textrm{Mult}(\widetilde H)/\Ker \varphi p \to \Htilde \varphi p \epi H$ is a non-split
central extension of $H$ such that $\Htilde \varphi p \subseteq G$.
The hypothesis on $G$ implies $\Htilde \varphi p=H$,
hence $\Ker \varphi p= \Mult(\emph{H})$ and $\Mult(\emph{H})\varphi \subseteq \Mult(\emph{G})$.
We can then apply Proposition \ref{mult-mult} and the theorem follows.
\Fp

As in \cite{RSV06}, if we assume that $H$ has no non-trivial central extensions then the following holds:
\begin{corollary}
\label{localization-superperfect}
Let $i:H\mono G$ be an inclusion of two non-abelian simple groups
and assume that $H$ is superperfect (i.e. $\textrm{\emph{Mult}}({H})=1$).
Let also $j:H=\Htilde  \mono \Gtilde $ denote
the induced homomorphism on the universal central extensions.
Then $i: H\mono G$ is a localization, if and only if $j:H\mono \Gtilde $ is a localization.  \Fp \bigskip
\end{corollary}

%%%%%%%%%%%%%%%%%%%%%%%%%%%%%%%%%%%%%%%%%%%%%%%%%%%%%%%%%%%%%%%%%
\section{Centers and localizations of quasi-simple groups}
%%%%%%%%%%%%%%%%%%%%%%%%%%%%%%%%%%%%%%%%%%%%%%%%%%%%%%%%%%%%%%%%%

Using ideas from Ol'shanskii \cite{Ol91} in \cite[Theorem B]{Obra94} (or \cite{Obra96}) Obraztsov showed the existence of simple complete groups
with a very special lattice of subgroups.  We will say that a group $G$ satisfies the \emph{Ol'shanskii's subgroup property},
(with respect to some family $\{G_i\mid i\in I\}$), if every proper subgroup of $G$ is either cyclic or contained in a conjugate of some $G_i$.
These strong results depend heavily on the assumption that $G$ has no involutions. For finite groups this is equivalent to saying that $G$ has odd order and thus must be solvable by the Feit--Thompson Theorem  \cite{FT63}, so the simple groups considered in this context must be infinite. We start with two important definitions (see \cite{Ol91}, Section 35):

\begin{definition}

Given a family $\{G_i\mid i\in I\}$ of groups, the \emph{free amalgam} of the groups $G_i$ is the set $\cup_{i\in I} G_i$, with $G_i\cap G_j=\{1\}$ for every $i\neq j$.

\end{definition}

Observe that no group structure is imposed, a priori, in a free amalgam.

\begin{definition}

Given a group $G$, a family $\{G_i\mid i\in I\}$ and the corresponding free amalgam $\mathcal{A}$, a map $\mathcal{A}\rightarrow G$ is an \emph{embedding} if it is injective and the restriction to every $G_i$ is a group homomorphism.

\end{definition}

In the sequel we will build on the mentioned results of Ol'shanskii and Obraztsov in order to generalize this notion of embedding to free products. With the notation of the previous definition, consider an embedding $f:\mathcal{A}\hookrightarrow G$, and the unique map $g:\mathcal{A}\rightarrow\ast_{i \in I}G_i$ that restricts to the inclusion $G_i\hookrightarrow \ast_{\in I}G_i$ for every $i\in I$. By the universal property of the free product of groups, there is a unique homomorphism $h:\ast_{\in I}G_i\rightarrow G$ such that $g h=f$. This fact motivates the following definition:

\begin{definition}

Let $F=\ast_{i\in I} G_i$ be the free product of the given family of groups. We call a homomorphism $\s: F\arr G$ a \emph{$*$-embedding}, if $\s\restr G_i: G_i\arr G$ is a monomorphism for each $i\in I$ and $G_i\s \cap G_j\s=1$ for all $i\ne j\in I$. This is equivalent to saying that $$\s\restr (\dot\bigcup_{i\in I}(G_i\setminus \{1\})\cup \{1\} \text{ is an injection. }$$

\end{definition}

By the previous observation, every embedding gives rise to a unique $*$-embedding. Observe however that a $*$-embedding is not necessarily injective, and in fact they can have very large kernels.

Before we present our first result concerning $*$-embeddings, we need one more definition. Recall that a group $G$ is {\it aspherical} (or {\it diagrammatically aspherical}) if there exists no reduced spherical diagram over $G$ consisting only of $0$-cells; a thorough introduction to the theory of diagrams is Chapter 4 of \cite{Ol91}.

Let us introduce now our first embedding result:

\begin{theorem}
\label{Obra}
Let $\{G_i \mid i\in I\}$, $(|I|\ge 2)$, be a countable family of countable groups
without involutions, of exponent $m > 10^{75}$ or
$m=\infty$. Then there is a $*$-embedding $F=\star_{i\in I} G_i\arr G$ such that the following holds.
\begin{enumerate}
\item[(1)] $G$ is simple, complete, aspherical and has no involutions;
\item[(2)] $G$ satisfies Ol'shanskii's subgroup property for $\{G_i \mid i\in I\}$;
\item [(3)]The Schur multiplier $\textrm{\emph{Mult}}({G})$ is a  free abelian group of infinite (countable) rank;
\item[(4)] If $m=\infty$, then $G$ is torsion-free and if $m<\infty$, then $G$ has exponent $m$.
\end{enumerate}
\end{theorem}
\Pf Let $\mathcal{A}$ be the free amalgam associated to the family $\{G_i \mid i\in I\}$, $(|I|\ge 2)$. According to above, to find a $*$-embedding $F=\star_{i\in I} G_i\arr G$ it is enough to find an embedding $\mathcal{A}=\star_{i\in I} G_i\arr G$ and then to extend it by linearity.

Now items (1), (2) and (4) are a direct consequence of \cite[Theorem B]{Obra94} (see also \cite[p. 386, Theorem 35.1]{Ol91}). For (3) note that the group $G$ of the theorem is aspherical, and hence has free abelian Schur multiplier of rank the size of the group $G$,  see \cite[p. 334, Theorem 31.1]{Ol91}.\\
\Fp

 \bigskip

We will derive an immediate

\begin{corollary}\label{Ismall} Assuming the statement of Theorem \ref{Obra},
with $\{G_i \mid i\in I\}$ $(\size {I}\ge 2)$ such that there are
no embeddings $G\mono G_i$ for any $i\in I$,  then  $G$  is co-hopfian. The additional hypothesis is satisfied if $G_i$ is hyperabelian or co-hopfian for $i\in I$.
\end{corollary}

\Pf Suppose for contradiction that $\va: G \mono G$ is a proper monomorphism, then $G\va \subsetneq G$, and by Ol'shanskii's subgroup property for $\{G_i \mid i\in I\}$ there are $i\in I$ and $g\in G$ such that $G\va\subseteq G_i^g$. Hence $\va (g^{-1})^*$ is an embedding of $G$ into $G_i$, which contradicts our assumption. Hence $G$ must be co-hopfian.

Assuming the additional hypothesis, we must show that an embedding $\psi: G\mono G_i$ does not exist. If $G_i$ is hyperabelian, then also $G$ must be hyperabelian and cannot be simple, as shown in Theorem \ref{Obra}(1). Otherwise $G_i$ must be co-hopfian and $G_i\subset G$ together with $\psi$ leads to a proper embedding $G_i\subset G\mono G\psi\subseteq G_i$, and so $G_i$ cannot be co-hopfian. Hence $\psi$ is not possible. \Fp \bigskip

We will also need a version of the above for $\size{I}\le 1$. If $I=\emptyset$, we will   replace $I$ by $I'=\{0,1\}$ and let $G_0=G_1=Z_m$ (with $m$ as in the theorem and $Z_m=Z$, if $m=\infty$). If $I=\{0\}$, then choose $G_0=H$ co-hopfian and $G_1=Z_m$. Thus we can apply Theorem \ref{Obra} and Corollary \ref{Ismall} for $I'$ (because $Z_m$ is hyperabelian and $H$ is co-hopfian) and get for any $I$ (also for $\size{I}= 0$ and $\size{I}= 1$) the following immediate

\begin{corollary}\label{cohI}
Let $H=1$ or $H$ be a countable co-hopfian group without involutions of any exponent $m > 10^{75}$ or
$m=\infty$. Then there is a group $G$ such that the following holds.
\begin{enumerate}
\item[(1)] $H\subseteq G$, with $G$ countable, simple, complete, aspherical and without involutions.
\item[(2)] If $H\ne 1$, then $G$ satisfies Ol'shanskii's subgroup property for $H$. If $H=1$, then $G$ satisfies Ol'shanskii's subgroup property for the empty set of groups.
\item [(3)]The Schur multiplier $\textrm{\emph{Mult}}({G})$ is a  free abelian group of infinite (countable) rank.
\item[(4)] If $m=\infty$, then $G$ is torsion-free, and if $m<\infty$, then $G$ has exponent $m$.
\item[(5)] $G$ is co-hopfian.
\end{enumerate}
\end{corollary}

\begin{theorem}
\label{Obra-localization}
Let $H$ be any simple, co-hopfian and complete group, which embeds properly into a simple, complete group $G$ satisfying Ol'shanskii's subgroup property with respect to $H$. Then $G$ is co-hopfian, and $H\subseteq G$ is a localization.
\end{theorem}
\Pf The group $G$ is co-hopfian by  Corollary \ref{Ismall}.  We can now use Corollary \ref{Criterion-complete} to show that $H\subseteq G$ is a localization.

Let $H'$ be any subgroup of $G$ isomorphic to $H$. If $G\cong H$, then $H\subset G\cong H$ contradicts co-hopfian. Hence $G\ncong H'$ and $H'\subsetneq G$ is not cyclic. By Ol'shanskii's subgroup property  there is $g\in G$ such that $H'\subseteq H^g$. Hence there is the embedding $\va(g^{-1})^*$ of $H'\cong H$ into $H$, which forces $H'=H^g$, because $H$ is co-hopfian. Thus (1) of Corollary \ref{Criterion-complete} follows.

Suppose that the centralizer $\fc_GH\neq 1$. If $H \fc_G H = G$, then  $H$ is normal in $G$, but $G$ is simple, which is impossible. Thus $H \fc_G H \ne G$ is not cyclic, and by Ol'shanskii's subgroup property, $H \fc_GH\subseteq H^g$ for some $g\in G$. Since $H$ has trivial center, it follows from $\fc_GH\neq 1$ that   $H\subsetneq H \fc_GH$ leads to a proper embedding of $H$ into itself (induced by $g$), which contradicts our assumption that $H$ is co-hopfian. Thus $\fc_GH=1$, part (2) of Corollary \ref{Criterion-complete} holds, and the corollary can be applied. Thus Theorem  \ref{Obra-localization} follows.  \Fp \bigskip

In the main theorem we shall need the following basic fact on central extensions. (Note that this holds for arbitrary cellular covers, where $\pi$ is assumed to be only an epimorphism.)
\begin{lemma}
\label{uniqueness} Suppose that $\Htilde \stackrel{\pi}{\epi} \Hbar \stackrel{\rho}{\epi} H$ and $\Gtilde \stackrel{\pi'}{\epi} \Gbar \stackrel{\rho'}{\epi} G$ are sequences of central extensions, where $\Htilde$ and $\Gtilde$ are the universal central extensions. Let $\psi, \psi': \Hbar \to \Gbar$ two homomorphisms.
\begin{itemize}
\item[(a)] If $\pi \psi=\widetilde\varphi \pi' = \pi\psi'$ for some $\widetilde\varphi: \Htilde\to \Gtilde$, then $\psi=\psi'$.
\item[(b)] If $\psi$ and $\psi'$ extend $\varphi:H\to G$, then $\psi=\psi'$.
\end{itemize}
\end{lemma}
 \Pf
Part (a) follows, since $\pi$ is epimorphism.
For (b) let $\widetilde\varphi: \Htilde \to \Gtilde$ be the unique extension of $\varphi$. Then also $\widetilde\varphi$ is the extension of both $\psi$ and $\psi'$, so the upper square in the following diagram commutes:
$$\xymatrix{
\Htilde \ar@{->>}[d]^{\pi} \ar@{.>}[r]^{\widetilde\varphi}  &             \Gtilde \ar@{->>}[d]^{\pi'} \\
\Hbar \ar@{->>}[d]^{\rho}  \ar@/^/[r]^\psi  \ar@/_/[r]^{\psi'}    &              \Gbar \ar@{->>}[d]^{\rho'} \\
H           \ar@{->}[r]^\varphi & G,
}
$$
%Finally $\pi\psi=\widetilde\varphi\pi'=\pi\psi'$  yields $\psi=\psi'$, because $\pi$ is an epimorphism.
Hence (a) applies, and we get $\psi=\psi'$.
\Fp

The next main result was already discussed in the Introduction.

\begin{theorem}
\label{main-theorem} Let $F$ be a free abelian group of countable rank and $A$ be any epimorphic image. Moreover, let $H$ be a countable, complete, co-hopfian simple group without involutions, and suppose that $m> 10^{75}$ or $m=\infty$ is such that all proper subgroups of $H$ are cyclic of an order dividing $m$, if $m<\infty$ and of infinite order otherwise.
\begin{itemize}
\item[(a)] Then there exists a localization $H\subset G$, with $G$ complete, simple, co-hopfian and Schur multiplier $\textrm{\emph{Mult}}({G})\cong F$.
\item[(b)] This localization lifts to a localization $\Htilde \to \Gtilde$ between the corresponding universal central extensions, where $\fz(\Gtilde)\cong F$.
\item[(c)] This localization lifts to another localization $\Hbar \subset \Gbar$, between the corresponding central extensions, where $\fz{\overline G}\cong A$.
\end{itemize}
\end{theorem}
\Pf
By Theorem \ref{Obra} there is an inclusion $H \subset G$ into a complete simple group $G$ of exponent $m$ with $\Mult(\emph{G})$ isomorphic to $F$. By Theorem \ref{Obra-localization} the group $G$ is co-hopfian, and this embedding is a localization, hence (a) holds.

Part (b) follows  directly from Theorem \ref{non-split}.
Suppose that $G$ contains a non-split central extension $\Hbar$ of $H$.
Then $\Hbar$ would be a non-cyclic proper subgroup of $G$,
which by the Ol'shanskii's subgroup property satisfies $\Hbar \subseteq H^g \cong H$
for some $g\in G$. Since $\Hbar$ is not simple, it must be a proper subgroup of $H$,  hence cyclic by the assumption on $H$, which is a contradiction. Therefore,
the additional hypothesis of Theorem \ref{non-split} holds, and we get that $\Htilde \to \Gtilde$ is a localization.

For part (c) we first define $\Gbar$ and $\Hbar$.
The epimorphism $\eta: F=\Mult(\emph{G})\epi \emph{A}$ induces a push-out diagram of central extensions:
$$\xymatrix{
0 \ar[r]  & F \ar@{->>}[d]^\eta \ar[r] & \Gtilde \ar@{->>}[d]^\pi \ar[r] & G \ar@{=}[d] \ar[r] &1\\
0 \ar[r] &  A \ar[r]           & \Gbar      \ar[r]^{p}         & G \ar[r] &1.}
$$
First note that $\pi^{-1} \fz\overline G$ is a normal subgroup of $\tilde G$ which contains $F$ and also $\pi^{-1} \fz\overline G\subseteq \Mult(\emph{G})=\emph{F}$, hence  $\fz\overline G=F\pi=A$.

Let $\Hbar:= (\Htilde) j \pi$, then we have a diagram of central extensions
$$\xymatrix{
\Htilde \ar@{->>}[d]^{\pi} \ar@{->}[r]^j  &             \Gtilde \ar@{->>}[d]^{\pi} \\
\Hbar \ar@{->>}[d]^{\rho}   \ar@{^(->}[r]^k    &              \Gbar \ar@{->>}[d]^{\rho} \\
H           \ar@{^(->}[r]^i & G,
}
$$
where $\Htilde \epi \Hbar$ is the universal central extension, and $k:\Hbar\subset \Gbar$ is the inclusion extending $i:H\subset G$. Note that $k$ is a strict inclusion $\Hbar \subsetneq \Gbar$, otherwise we would have $\Hbar\cong \Gbar$, and thus $H\cong \Hbar/\fz(\Hbar)\cong \Gbar/\fz(\Gbar)\cong G$. This would yield a strict monomorphism $H\mono G\cong H$, contradicting  that $H$ is co-hopfian.

For (c) it remains to show that $\Hbar\subset \Gbar$ is the desired localization.
Thus let $\vabar: \Hbar \to \Gbar$ be any homomorphism.

If $\vabar=0$, then $0\in \End \Gbar$ extends $\vabar$, and also $\widetilde\varphi=0: \Htilde \to \Gtilde$ extends uniquely $\vabar=0$. If $\psibar: \Gbar \to \Gbar$ is a homomorphism which extends $\vabar=0$, and $\widetilde \psi:\Gtilde \to \Gtilde$ is the unique extension of $\psibar$, then $j \widetilde \psi =\widetilde \varphi=0$ by the uniqueness property of universal central extensions.
$$
\xymatrixcolsep{5pc}\xymatrix@!0{
 \widetilde H   \ar@{->}[dr]_{\widetilde \varphi=0}  \ar@{^{(}->}[rr]^j \ar@{->>}[dd]^\pi & & \widetilde G\ar@{->}[dl]_{0}  \ar@/^/[dl]^{\widetilde \psi} \ar@{->>}[dd]^(.6)\pi
  \\
  &  \widetilde G    \ar@{->>}[dd]^(.6)\pi
\\
\bar H   \ar@{->}[dr]_{\bar \varphi=0}  \ar@{^{(}->}'[r][rr]^(.3)k && \bar G  \ar@{->}[dl]_{0}  \ar@/^/[dl]^{\bar \psi} &
\\
  &  \bar G
}
$$
The the uniqueness property for the localization $j$ we have $\widetilde\psi=0$, and therefore $\pi\bar\psi= \pi 0$. Since $\pi$ is an epimorphism we obtain $\bar \psi=0$ as desired (see Lemma \ref{uniqueness} (a)).

Suppose finally that $\vabar \neq 0$. Hence $\Ker \vabar \subseteq \fz(\Hbar)$, because $\Hbar/\fz(\Hbar)\cong H$ is simple.
Then $\Ker\vabar\rho\neq \Hbar$, otherwise $\vabar\rho= 0 \rho$ and this would imply that both $\vabar$ and $0$ extend the trivial map $0: H\to G$, and then $\vabar=0$, by Lemma \ref{uniqueness} (b).
We have that $\Ker \vabar \rho \subseteq \fz(\Hbar)$,  because $H$ is simple.

Let $S=\Hbar\vabar\rho\cong \Hbar/\Ker \vabar \rho\neq 0$ and consider the canonical epimorphism $S\epi H$, induced by the inclusion $\Ker \vabar \rho \subseteq \fz(\Hbar)$.
In particular, $S$ cannot be cyclic. By Ol'shanskii's subgroup property it follows that either $S=G$ or $S\subseteq H^g\subset G$ for some $g\in G$.

In the first case we have $S=G$, thus an epimorphism $G=S\epi H$, which must be an isomorphism, as $G$ is simple. This yields a proper inclusion $H\varsubsetneq G\cong H$ of $H$, which contradicts the assumption that $H$ is co-hopfian.

 In the second case  we have a monomorphism $S\subseteq H^g\cong H$ and the epimorphism $S\epi H$, as before. But, by hypothesis,  proper subgroups of $H$ are cyclic, thus $S=H^g$.

 If $\a: H\arr S$ is the given isomorphism,
 then $H\a= H^g$ and $\a(g^{-1})^*\in \Aut H$.
 Since $H$ is complete ($\Aut H=H^*$), we find $h\in H$ such that $\a(g^{-1})^*=h^*$, hence $\a= (hg)^*$.
 Denote $x=hg$ and choose preimages $y\in \Hbar$, and $z\in \Htilde$ such that $y\rho=x$, and $z \pi =y$.

$$
\xymatrixcolsep{5pc}\xymatrix@!0{
& \Htilde \ar@{->}[dl]_{z*}^(.3)\cong \ar@{->}[rr]^j \ar@{->>}'[d][dd]^\pi \ar@{->}[dr]^{\tilde \varphi} & & \Gtilde \ar@{->>}[dd]^\pi \ar@{->}[dl]_{z^*} \ar@/^/[dl]^{\widetilde\psi}
\\
\widetilde H^z \ar@{->}[rr] \ar@{->}[dd] &  & \Gtilde \ar@{->>}[dd]
\\
& \bar H  \ar@{->>}[dl]_(.3){y^*}^(.3)\cong \ar@{->}[dr]^{\bar \varphi}  \ar@{^{(}->}'[r][rr]^k \ar@{->>}'[d][dd]^\rho & & \bar G   \ar@{->>}[dd]^\rho \ar@{->}[dl]_(.6){y^*} \ar@/^/[dl]^{\bar\psi}
\\
\bar H^y \ar@{^{(}->}[rr] \ar@{->>}[dd] &  &  \bar G    \ar@{->>}[dd]^(.6)\rho
\\
&H    \ar@{->}[dl]_(.3){x^*}^(.3)\cong \ar@{^{(}->}'[r][rr]^(.3)i \ar@{->}[dr]^{\varphi} && G  \ar@{->}[dl]^\cong _{x^*} &
\\
S=H^g  \ar@{^{(}->}[rr]  &&  G
}
$$
Now we can apply Lemma \ref{uniqueness}(b) to infer that $\bar \varphi= y^* \restr \Hbar$.
Thus clearly,  $y^*: \Gbar \to \Gbar$ extends $\vabar$, i.e. $k y^*=\vabar$.

We show next that this extension is unique.
If $\psibar:\Gbar \to \Gbar$ is another homomorphism such that  $k \psibar = \vabar$, then $\psibar$ lifts to a unique $\widetilde \psi$ such that $\widetilde\psi \pi = \pi \psibar$. Now $j z^*= \widetilde \varphi$ and also $j \widetilde\psi=\widetilde\varphi$. For the second equality use that $j\widetilde\psi=j\pi\bar \psi= k\pi \bar\psi=\pi\varphi$, and the uniqueness of liftings to universal central extensions. Hence, the uniqueness of extensions of the localization $j$ yields $\widetilde\psi =z^*$. Thus we have $\pi \bar\psi = \pi y^*$, and since $\pi$ is an epimorphism, we get $\psibar=y^*$. Thus $\Hbar\subset\Gbar$ is also a localization, as desired. \Fp \bigskip

Assuming that there are such groups $H$ with trivial Schur multiplier, we get the following corollary.

\begin{corollary}
Let $A$ be any countable abelian group and let $H$ be any countable, non-abelian simple group without involutions of exponent  $m> 10^{75}$ or $m=\infty$.
Assume that $H$ is co-hopfian, complete, with Schur multiplier $\textrm{\emph{Mult}}({H})=1$.
Then there exists a localization $H \subset \Gbar$  with $\fz \Gbar=A$, and $G$ is a complete, simple and co-hopfian group of exponent $m$. \Fp \bigskip
\end{corollary}

\Pf From Corollary \ref{cohI} we obtain an inclusion $H\subset G$ which is a localization by Theorem \ref{Obra-localization} (or Theorem \ref{main-theorem}(b)). Corollary \ref{trivial-schur} says that $H\subset \Gtilde$ is also a localization.
By the proof of part (c) of Theorem \ref{main-theorem} (when $\Htilde=\Hbar=H$), we notice that it is not necessary to assume that all proper subgroups of $H$ are cyclic of exponent $m$. Indeed, following the argument of the proof, the epimorphism $S\epi H$ is an isomorphism, and $S$ cannot be a proper subgroup of $H^g$, because otherwise we would have a strict monomorphism $H\cong S \subsetneq H^g\cong H$, which would contradict that $H$ is co-hopfian. The proof can then continue without this additional assumption. Hence we get a localization $H=\Hbar \subset \Gbar$, with center $\fz \Gbar=A$,   satisfying the desired properties.  \Fp \bigskip

In view of the last corollary it would be interesting to answer the following question.

\begin{problem}
\label{Ol}
Is there a simple, complete co-hopfian group $H$ without involutions and   trivial Schur multiplier?\end{problem}

\begin{remark}
{\rm
In Theorem \ref{main-theorem} we can replace the property that all proper subgroups of  $H$ are cyclic by the more general property that $H$ does not contain any non-split central extension of itself as a proper subgroup.

According to Proposition \ref{universalcover-localization} the simple groups satisfying this property are exactly those groups $H$ whose universal central extension $\Htilde \epi H$ is a localization.

In the theory of the localization of groups (and other categories) we sometimes find morphisms that are localizations and cellular covers (or co-localizations) at the same time; see e.g. \cite{FGSS07} and \cite{GRS12}.
}\end{remark}

%%%%%%%%%%%%%%%%%%%%%%%%%%%%%%%%%%%%%%%%%%%%%%%%%%%%%%%%%%%%%%%%%
%                       BIBLIOGRAPHY                            %
%%%%%%%%%%%%%%%%%%%%%%%%%%%%%%%%%%%%%%%%%%%%%%%%%%%%%%%%%%%%%%%%%
%\textbf{Acknowledgments.} The second-named author wish to thank the hospitality of Prof. R. G\"{o}bel and his help in the elaboration of most of this work.

\vskip 0.5 cm

\setlength{\baselineskip}{0.6cm}

Ram\'{o}n Flores

Departamento de Geometr\'{i}a y Topolog\'{i}a, Universidad de Sevilla-IMUS

E-41012 Sevilla, Spain, e-mail: {\tt ramonjflores@us.es}

\bigskip\noindent

Jos\'e L. Rodr\'{\i}guez

Departamento de Matem\'aticas, Universidad de Amer\'{\i}a

E-04120 La Ca\~nada de San Urbano,  Almer\'{\i}a,
 Spain, e-mail: {\tt jlrodri@ual.es}

\end{document}